\documentclass[a4paper,reqno]{amsart}
\usepackage{amsmath,amssymb,graphicx,lscape}
\newcommand{\KH}{{\it KH}}
\newcommand{\Kh}{{\it Kh}}
\newcommand{\ta}{\tilde{a}}
\newcommand{\tb}{\tilde{b}}
\newcommand{\tD}{\widetilde{D}}
\newcommand{\tn}{\tilde{n}}
\newcommand{\bdots}{\reflectbox{$\ddots$}}%graphicx
\newcommand{\bs}{\reflectbox{$/$}}
\DeclareMathOperator{\supp}{Supp}

\newtheorem{thm}{Theorem}[section]
\newtheorem{prop}[thm]{Proposition}
\newtheorem{lem}[thm]{Lemma}

\theoremstyle{definition}
\newtheorem{define}[thm]{Definition}
\theoremstyle{remark}
\newtheorem{rem}[thm]{Remark}
\newtheorem{ex}[thm]{Example}

\title{Khovanov homology and Rasmussen's $s$-invariants for pretzel knots}
\author{RYOHEI SUZUKI}

\begin{document}

\begin{abstract}
We calculated the rational Khovanov homology of some class of pretzel knots, %
by using the spectral sequence constructed by P. Turner. %
Moreover, we determined the Rasmussen's $s$-invariant of almost of pretzel %
knots with three pretzels.
\end{abstract}

\maketitle

\section{Introduction and main results}

In \cite{Kho}, Khovanov constructed a bi-graded homology invariant $\KH(L)$ of links, %
whose graded Euler characteristic is the unnormalized Jones polynomial of $L$. %
Lee \cite{Lee} modified this theory and constructed another invariant $\KH_{Lee}(L)$, %
which is a singly-graded homology. %
Lee showed that the rational Khovanov homology is representable as the $E_2$ term of a spectral sequence %
which converges to $\mathbb{Q}^{{}\oplus 2^n}$, where $n$ is a number of components of a link $L$. %
In \cite{Ras}, Rasmussen defined a knot invariant $s(K)$ using the Khovanov homology and Lee theory. %
This invariant $s$ gives a lower bound for the slice genus. %
There are several inequalities for $s(K)$: %
slice-Bennequin inequality proved by Shumakovitch \cite{Shu} and Plamenevskaya \cite{Pla} independently, %
sharper slice-Bennequin inequality by Kawamura \cite{Kawamura}, %
and crossing change inequality by Livingston \cite{Liv} and Rasmussen \cite{Ras}. %

There was no theoretical computational tool for Khovanov homology %
but the \textit{skein exact sequence} obtained from a short exact sequence of Khovanov homology
\begin{equation}
0 \to C(D'') \to C(D) \to C(D') \to 0
\end{equation}
where $D$ is a link diagram and $D'$ and $D''$ are obtained by resolving the same crossing of $D$ %
to the 0-smoothing and the 1-smoothing, respectively. %
We would have used this long exact sequence repeatedly before, but it requires careful bookkeeping. %
Turner \cite{PT} defined a spectral sequence, which is obtained by tying up in a bundle the skein exact sequences. %
It is simple enough to enable us fast theorical computation. %
In \cite{PT}, Turner gives an example of an application to $(3,n)$-torus links. %

In this paper, we give another application of Turner's spectral sequence to pretzel knots. %
Our main results are stated below. %
Theorem \ref{Thm:main} is a result for $\KH(K)$, %
and Theorem \ref{Thm:main2} and \ref{Thm:main3} are for $s(K)$. %
In Section 2, we summarize without proofs some basic facts. %
In Section 3, we prove Theorem \ref{Thm:main}. %
In Section 4, we give remarks about the assumption of Theorem \ref{Thm:main}. %
In Section 5, we prove Theorem \ref{Thm:main2} and \ref{Thm:main3}. %

\begin{thm}\label{Thm:main}
Suppose $p$ is an odd number and $p \ge 9$, $q=p-2$, %
and $r$ is a positive even number. %
Let $K$ be the pretzel knot $P(p,-q,-r)$. %
Then Rasmussen's $s$-invariant $s(K)$ is given by
\begin{equation}
s(K)=2.
\end{equation}
Furthermore, the rational Khovanov homology $\KH(K)$ is computable %
and given explicitly as in formulae (\ref{eq:000}), (\ref{eq:001}), (\ref{eq:002}).
\end{thm}

In the proof of Theorem \ref{Thm:main}, %
we obtain the value of $s$-invariant before we obtain the homology. %
Hence our approach is easily applied to calculation of $s$-invariant of more larger class of pretzel knots. %
In fact, by combinatorial use of sharper slice-Bennequin inequality and this approach, %
we could obtain the value of $s$-invariant of almost all of pretzel knots with three pretzels. %

\begin{thm}\label{Thm:main2}
Suppose $p$ and $q$ are odd numbers and $p \ge 3 $, $q \ge 3$, %
and $r$ is a positive even number. %
Let $K$ be the pretzel knot $P(p,-q,-r)$. %
Then Rasmussen's $s$-invariant $s(K)$ is given by
\begin{equation}
s(K)=p-q.
\end{equation}
\end{thm}

\begin{thm}\label{Thm:main3}
Suppose $p$, $q$ and $r$ are odd numbers and $p \ge 3 $, $q \ge 3$, $r \ge 3$. %
Let $K$ be the pretzel knot $P(p,-q,-r)$. %
Then Rasmussen's $s$-invariant $s(K)$ is given by
\begin{equation}
s(K)=%
\begin{cases}
0 & \text{if }  p > \min\{q,r\} \\
2 & \text{if } p < \min\{q,r\} \\
% ??? & \text{if } p = \min\{q,r\}
\end{cases}
\end{equation}
\end{thm}

We note some properties of pretzel knot with three pretzels and its $s$-invariant. %

\begin{prop}\label{prop:property}
i) A pretzel knot $P(p,q,r)$ is indeed knot %
if and only if all of $p$, $q$, $r$ are odd numbers, or %
one of these integers is even number and the other two are odd numbers. %

ii) For all $\sigma \in {\mathfrak{S}}_3$, $s_1,s_2,s_3 \in \mathbb{Z}$, %
\begin{equation}
 P(s_1,s_2,s_3) \simeq P(s_{\sigma(1)},s_{\sigma(2)},s_{\sigma(3)})
\end{equation}
and %
\begin{equation}
 s(P(-s_1,-s_2,-s_3)) = s(P(s_1,s_2,s_3)^!) = -s(P(s_1,s_2,s_3)).
\end{equation}
\end{prop}

By Proposition \ref{prop:property}, we do not need to consider all patterns. %
For example, if $p$, $q$ and $r$ satisfy the assumptions of Theorem \ref{Thm:main2}, then $s(P(p,-q,r))=p-q$.

\section{Tools and notation}
\subsection{Turner's spectral sequence}

Let $L$ be an oriented link and $D$ be an oriented diagram of $L$. %
Let $ \{ 1, \dots , m \} $ be a subset of the set of the crossings of $D$. %
For $1 \le k \le m$, let $D_{(k)}$ be the diagram obtained from $D$ by resolving the crossing $ 1,\dots,k$ to 1-smoothings, %
and $\tD_{(k)}$ the diagram obtained from $D$ %
by resolving the crossing $ 1,\dots,k-1$ to 1-smoothings and the crossing $k$ to a 0-crossing. %
For $k=0$, let $ D_{(0)} = \tD_{(0)} = D $ as an oriented diagram. 

If the $k$th crossing is positive then there is the canonical orientation of the diagram $\tD_{(k)}$: %
since the 0-smoothing is the oriented resolution, $\tD_{(k)}$ inherits an orientation from $D_{(k-1)}$. %
But for $D_{(k)}$, there is no canonical one, so we can choose any orientation. %
On the contrary, if the $k$th crossing is negative then $D_{(k)}$ inherits an orientation from $D_{(k-1)}$ %
and we can choose any orientation for $\tD_{(k)}$. 

Turner's spectral sequence, which is converging to $\KH(D)$, %
is obtained by arranging $ \KH(\tD_{(k)}) $'s and $ \KH(D_{(m)}) $ at the appropriate position %
represented by some constants $\ta_k$, $\tb_k$, $A_k$, and $B_k$ defined from $D_{(k)}$'s and $\tD_{(k)}$'s.

\begin{prop}[P. Turner \cite{PT} Proposition 2.2]\label{Prop:Tur1}
For every fixed $j \in \mathbb{Z}$, %
there is a spectral sequence $ ( E_r^{*,*},d_r : E_r^{s,t} \to E_r^{s+r,t-r+1} ) $ %
converging to $ \KH^{*}_{j}(D) $ with $ E_1 $-page given by
\begin{equation}
E_1^{s,t} =
\begin{cases}
 \KH^{s+t+A_s+\ta_{s+1}}_{j+B_s+\tb_{s+1}}(\tD_{(s+1)}) & s=0, \dots, m-1\\
 \KH^{m+t+A_m}_{j+B_m}(D_{(m)}) & s=m \\
 0 & \text{otherwise}
\end{cases}
\end{equation}
where $\ta_k$, $\tb_k$, $A_k$, and $B_k$ are the constants defined by the following:
\begin{align*}
n^{+}_k &= \text{number of positive crossing in } D_{(k)}\\
n^{-}_k &= \text{number of negative crossing in } D_{(k)}\\
\tn^{+}_k &= \text{number of positive crossing in } \tD_{(k)}\\
\tn^{-}_k &= \text{number of negative crossing in } \tD_{(k)}\\
\end{align*}
\begin{align*}
a_k &= n^{-}_{k-1} - n^{-}_{k} -1 &\text{and}& &b_k &= 3 a_k + 1 \\
\ta_k &= \tn^{-}_{k-1} - \tn^{-}_{k} &\text{and}& &\tb_k &= 3 \ta_k -1.
\end{align*}
\begin{equation*}
A_0 = B_0 = 0, \, A_k = \sum_{i=1}^{k} a_i , \text{ and } B_k = \sum_{i=1}^{k} = 3 A_k + k.
\end{equation*}
\end{prop}

\subsection{Lee theory}

To prove the theorems, we need to make use of Lee theory. %
Lee theory is a variant of rational Khovanov homology %
obtained from the same underlying vector spaces %
but using a different differential. %
We denote it by $\KH_{Lee}^*(L)$. %
We summarize the results we need about Lee theory in the following proposition.

\begin{prop}\label{Prop:Lee}
For a knot $K$,

(i) $\KH_{Lee}^t(K) =
\begin{cases}
\mathbb{Q} \oplus \mathbb{Q} & t = 0 \\
0 & t \neq 0
\end{cases}$

(ii) There is a spectral sequence $(E_r^{*,*},d_r)$ converging to $\KH_{Lee}^{*}(K)$ with $E_1^{s,t} = \KH^{s+t}_{4s\pm1}(K)$. %
The differential $d_r$ has bi-degree $(1,4r)$ in term of index of the Khovanov homology.
\end{prop}

\subsection{Notation}

 Let $ P(p,q,r) $ be the standard diagram of $(p,q,r)$-pretzel link (or link itself). %
By abuse of notation, if $p$ (or $q$,$r$) is equal to 0, %
we consider that there is no half twist at the corresponding pretzel. %
For example, $ P(p,q,0) \simeq T(2,p) \# T(2,q) $ as an unoriented link. 

Let $ P(p,q,r^{\vee}_{\wedge}) $ be the diagram obtained %
from the diagram $ P(p,q,r+\epsilon) $ %
by resolving the crossing at the top of the corresponding pretzel %
to a $x$-smoothing, where $ \epsilon =1,x=1 $ if $r \ge 0$, %
$ \epsilon = -1 ,x=0$ if $ r<0 $. %
Obviously, $ P(p,q,r^{\vee}_{\wedge}) \simeq P(p,q,0) $ as an unoriented link. %
$ P(p^{\vee}_{\wedge},q,r),P(p,q^{\vee}_{\wedge},r) $ is also defined in the same manner.

$\KH(L)$ and $\Kh_{\mathbb{Q}}(L)$ denotes the rational Khovanov homology of $L$ and its Poincar\'{e} polynomial, respectively. 

\section{Proof of Theorem 1.1}

We prove the theorem \ref{Thm:main} by two steps. Each step consists of three small steps : %
first, we write down the spectral sequence, %
secondly, we describe the possibility of $\KH(K)$, %
thirdly, we determine $\KH(K)$ (by using Lee theory).

\begin{prop}[Step 1]\label{prop:step1}
If $p:\text{odd} \ge 9$, $q = p-2 $,
\begin{multline}
\Kh_{\mathbb{Q}}(P(p,-q,0)) =%
 q^{-2p+5}t^{-p+2} + q^{-2p+9}t^{-p+3} + 2q^{-2p+9}t^{-p+4} + q^{-2p+11}t^{-p+5} \\
 + \sum_{\substack{n=-p+7 \\ n:\text{even}}}^{-2}q^{2n}t^{n}((I_n+1)q^{-1}t^{-2} + (I_n+2)q^{-1}t^{-1} + I_n q^{1}t^{-1} + (I_n+1)q^{1}t^{0}) \\
 + (I_0+3)q^{-1}t^{-2} + (I_0+3)q^{-1}t^{-1} + (I_0+2)q^1t^{-1} + (I_0+4)q^1t^0 \\
 + (I_0+4)q^3t^0 + (I_0+3)q^3t^1 + (I_0+3)q^5t^1 + (I_0+4)q^5t^2 \\
 + \sum_{\substack{m=2 \\ m:\text{even}}}^{p-3}q^{2m}t^{m}((I_m+1)q^3t^0 + I_m q^3t^1 + (I_m+2)q^5t^1 + (I_m+1)q^5t^2) \\
 + q^{2p+3}t^p. \label{eq:000}
\end{multline}
where $ I_m = \frac{p-m-3}{2}, I_0 = \frac{p-9}{2}, I_n = \frac{p+n-5}{2} $.
In other words, $ \KH_\mathbb{Q}(P(p,-q,0)) $ is given by Table \ref{TabHuge1} (see page \pageref{TabHuge1}).
\end{prop}

\begin{proof}
We use Proposition \ref{Prop:Tur1} (in case $m=q-1$).
Let $D = D_{(0)} = \tD_{(0)} = P(p,-q,0)$ and its orientation be as shown in Figure \ref{fig:ori1}. %
For $1 \le k \le q-1$, $D_{(k)} := P(p,-(q-k),0)$ (which inherits the orientation from $D_{(k-1)}$), %
$\tD_{(k)} := P(p,-(q-k)^{\vee}_{\wedge},0)$ and its orientation be as shown in Figure \ref{fig:ori2}. %
It is easy to see that $\tD_{(k)} \simeq T(2,p)$ for $1 \le k \le q-1$, $D_{(q-1)} \simeq T(2,p)$, %
$\ta_s = q-s+1$, $\tb_s = 3q-3s+2$, $A_s = 0$, $B_s = s$ for $1 \le s \le q-1$.

\begin{figure}[htbp]
\begin{center}
\begin{minipage}{5cm}
\includegraphics{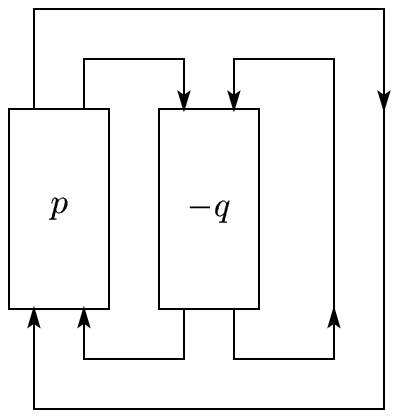}
\caption{}\label{fig:ori1}
%orientation of $P(p,-q,0)$
\end{minipage}
\begin{minipage}{5cm}
\includegraphics{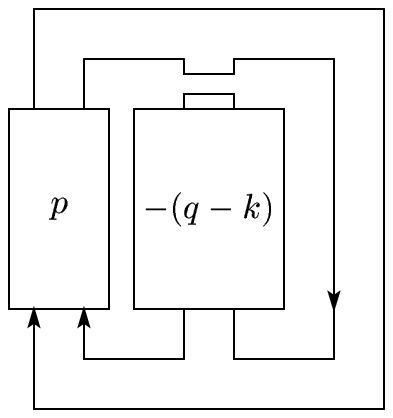}
\caption{}\label{fig:ori2}
%orientation of $P(p,-(q-k)^{\vee}_{\wedge},0)$
\end{minipage}
\end{center}
\end{figure}

The spectral sequence is
\begin{align}
E_1^{s,t} &= \KH^{t+q}_{j+3q-2s-1}(T(2,p)) \text{ if } s=0,\dots,q-2\\
E_1^{q-1,t} &= \KH^{t+q-1}_{j+q-1}(T(2,p)).
\end{align}
From the result of the (integral) Khovanov homology of $(2,p)$-torus link (see \cite{Kho}), %
one can write down the spectral sequence for each $j$. %
For $j \ge 2p+1$, $j=2p-3$ or $j \le 5-2p$, the spectral sequence is concentrated in one column %
and hence collapses for dimensional reasons. %
For other $j$, there exists one or several possible (nontrivial) $d_r$'s ($r\ge1$). %

If a generator of $E_1^{s,t}$ for some $j$ survives to $E_\infty$, %
it corresponds to a generator of $ \KH^{s+t}_{j}(D) $. %
So the generator of $E_1^{s,t}$ for some $j$ located at the initial/terminal end of these differential %
(we note that it may possibly survive to $E_\infty$) %
corresponds to a {\it possible} generator of $ \KH^{s+t}_{j}(D) $. 

Now the situation is as shown in Table \ref{TabHuge2} (see page \pageref{TabHuge2}). %
We note that there is no possible generator outside the two diagonals: $\deg = 2\dim+2 \pm 1$. %
We require the support of Lee theory. %
Suppose that the possible generator $z$ in bi-degree $(-p+2,-2p+7)$ is indeed a generator. %
Then it survives to $E_\infty$ of Lee's spectral sequence, (there is nothing which kills $z$), %
so $\KH_{Lee}^{-p+2} \neq 0$, which contradicts Proposition \ref{Prop:Lee}. Hence this generator $z$ must be {\it fake} (not a generator). %
Moreover, by looking back to Turner's spectral sequence, the possible generator $w$ in bi-degree $(-p+3,-2p+7)$ must be also fake. %
This is from the following reason. %
Since $z$ is fake, and the $E_1$ for $j=-2p+7$ is %
$\begin{array}{|c||c|c|} \hline
t\bs s & 0 & 1 \\ \hline \hline
-p+2& 1 & 1 \\ \hline
\end{array}$\,, %
the generator $\bar{z} \in E_1^{0,-p+2}$ (which corresponds $z$) must be killed at $E_1$, %
and when $\bar{z}$ is killed, another generator $\bar{w} \in E_1^{1,-p+2}$ must be killed at the same time, which corresponds to $w$. %
Everything is going on like this. %
By using two spectral sequences alternately, we can show that %
for each even number $3-p \le n \le -2$, exactly 1 pair of possible generators in bi-degree $(n-1 _{\text{ and }} n,1+2n)$ must be fake, %
and other possible generators must be indeed generators (details are omitted). %
Thus we obtain the result as shown in the proposition.
\end{proof}

\begin{rem}
As mentioned in the proof, %
there is no possible generator outside the two diagonals $\deg = 2\dim+2 \pm 1$. %
This indicates that the rational Khovanov homology of $P(p,-q,0)$ is computable if we just know the Jones polynomial of $P(p,-q,0)$, %
which is equal to $V(T(2,p))V(T(2,-q))$.
\end{rem}

\begin{prop}[Step 2 in case $r=2$]
If $p:\text{odd} \ge 9$, $q = p-2 $, $r=2$, $K=P(p,-q,-r)$ then
\begin{gather}
\Kh_{\mathbb{Q}}(K) = q^4 t^2 \Kh_{\mathbb{Q}}(P(p,-q,0)) + q^1 + q^3 + q^3 t^1 + q^9 t^3 \label{eq:001} \\
s(K) = 2
\end{gather}
and $K$ is H-thin.
\end{prop}

\begin{figure}[htbp]
\begin{center}
\begin{minipage}{5cm}
\includegraphics{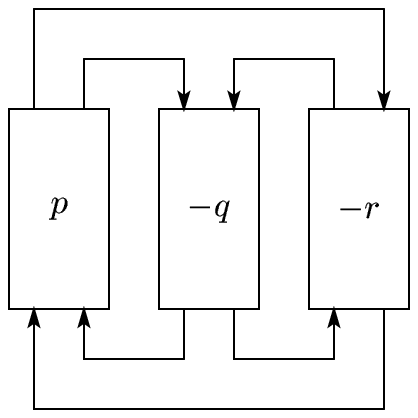}
\caption{}\label{fig:ori3}
%orientation of $P(p,-q,-r)$
\end{minipage}
\begin{minipage}{5cm}
\includegraphics{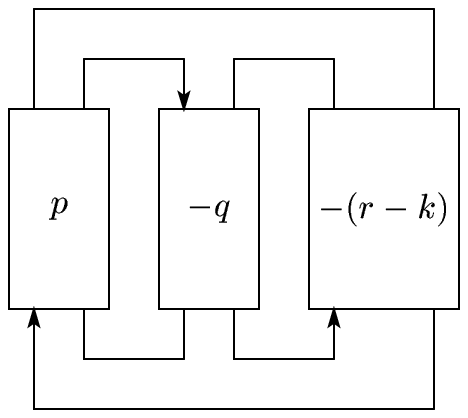}
\caption{}\label{fig:ori4}
%orientation of $P(p,-q,-(r-k))$
\end{minipage}
\end{center}
\end{figure}

\begin{proof}
We use Proposition \ref{Prop:Tur1} (in case $m=2$).
Let $D = P(p,-q,-r)$ and its orientation be as shown in Figure \ref{fig:ori3}, %
$D_{(k)} = P(p,-q,-(r-k)) $ for $k=1,2$ and its orientation as shown in Figure \ref{fig:ori4}. %
We note that the orientations of unmarked arcs vary with whether $k=1$ or $k=2$. %
$\tD_{(k)} = P(p,-q,-(r-k)^{\vee}_{\wedge}) $ (which inherits the orientation from $D_{(k-1)}$) for $k=1,2$. %
Then $ \ta_1 = \ta_2 = 0$, $\tb_1 = \tb_2 = -1$, $A_1 = -3$, $B_1 = -8$, $A_2 = -2$, $B_2 = -4 $.

It is easy to see $ \tD_{(1)} \simeq T(2,2) $ and $ \tD_{(2)} $ is equivalent to the link $L$ obtained from $T(2,2)$ %
by changing the orientation of one component. %
In this case, we may say that $L \simeq T(2,-2)$. %
We can easily get the homology of $L$ from the homology of $T(2,2)$ by shifting indexes: %
\begin{equation}
\KH^i_j(L) \cong \KH^{i+2}_{j+6}(T(2,2)).
\end{equation}

The spectral sequence is
\begin{align}
E_1^{0,t} &= \KH^{t}_{j-1}(T(2,2)) \\
E_1^{1,t} &= \KH^{t-2}_{j-9}(L) = \KH^{t}_{j-3}(T(2,2)) \\
E_1^{2,t} &= \KH^{t}_{j-4}(P(p,-q,0)).
\end{align}
Since the homology of $T(2,2)$ is very simple, the spectral sequence is also simple. %
This fact enables us to carry out the computation.
For $ j \neq 1,3,5,7,9 $, the spectral sequence is concentrated in one column $s=2$, %
hence collapses for dimensional reasons. So we obtain $ \KH^{i}_{j}(P(p,-q,-r)) \cong \KH^{i-2}_{j-4}(P(p,-q,0)) $. %
For $j=1,3,5,7,9$, the spectral sequence is as shown in Table \ref{Tab:tss1}. %
(Let $I_0:=(p-9)/2$.)

\begin{table}[ht]
\begin{center}
\begin{minipage}[t]{3.4cm}
\small
\begin{tabular}{cccc}
j=9 &  &  &  \\ \hline
\multicolumn{1}{|c||}{t\bs s} & \multicolumn{1}{c|}{0} & \multicolumn{1}{c|}{1} & \multicolumn{1}{c|}{2} \\ \hline\hline
\multicolumn{1}{|c||}{2} & \multicolumn{1}{c|}{} & \multicolumn{1}{c|}{$1$} & \multicolumn{1}{c|}{$I_0+4$} \\ \hline
\multicolumn{1}{|c||}{1} & \multicolumn{1}{c|}{} & \multicolumn{1}{c|}{} & \multicolumn{1}{c|}{$I_0+3$} \\ \hline
 &&& \\
 &&& \\
\end{tabular}%
\end{minipage}\quad
%%%%%%%%%%%%%%%%%%%%%%%%%%
\begin{minipage}[t]{3.4cm}
\small
\begin{tabular}{cccc}
j=7 &  &  &  \\ \hline
\multicolumn{1}{|c||}{t\bs s} & \multicolumn{1}{c|}{0} & \multicolumn{1}{c|}{1} & \multicolumn{1}{c|}{2} \\ \hline\hline
\multicolumn{1}{|c||}{2} & \multicolumn{1}{c|}{$1$} & \multicolumn{1}{c|}{$1$} & \multicolumn{1}{c|}{} \\ \hline
\multicolumn{1}{|c||}{1} & \multicolumn{1}{c|}{} & \multicolumn{1}{c|}{} & \multicolumn{1}{c|}{$I_0+3$} \\ \hline
\multicolumn{1}{|c||}{0} & \multicolumn{1}{c|}{} & \multicolumn{1}{c|}{} & \multicolumn{1}{c|}{$I_0+4$} \\ \hline
 &&& \\
\end{tabular}%
\end{minipage}\quad
%%%%%%%%%%%%%%%%%%%%%%%%%%
\begin{minipage}[t]{3.4cm}
\small
\begin{tabular}{cccc}
j=5 &  &  &  \\ \hline
\multicolumn{1}{|c||}{t\bs s} & \multicolumn{1}{c|}{0} & \multicolumn{1}{c|}{1} & \multicolumn{1}{c|}{2} \\ \hline\hline
\multicolumn{1}{|c||}{2} & \multicolumn{1}{c|}{$1$} & \multicolumn{1}{c|}{} & \multicolumn{1}{c|}{} \\ \hline
\multicolumn{1}{|c||}{1} & \multicolumn{1}{c|}{} & \multicolumn{1}{c|}{} & \multicolumn{1}{c|}{} \\ \hline
\multicolumn{1}{|c||}{0} & \multicolumn{1}{c|}{} & \multicolumn{1}{c|}{$1$} & \multicolumn{1}{c|}{$I_0+4$} \\ \hline
\multicolumn{1}{|c||}{\llap{$-$}1} & \multicolumn{1}{c|}{} & \multicolumn{1}{c|}{} & \multicolumn{1}{c|}{$I_0+2$} \\ \hline
\end{tabular}%
\end{minipage}\\
%%%%%%%%%%%%%%%%%%%%%%%%%%
\begin{minipage}[t]{3.4cm}
\small
\begin{tabular}{cccc}
j=3 &  &  &  \\ \hline
\multicolumn{1}{|c||}{t\bs s} & \multicolumn{1}{c|}{0} & \multicolumn{1}{c|}{1} & \multicolumn{1}{c|}{2} \\ \hline\hline
\multicolumn{1}{|c||}{0} & \multicolumn{1}{c|}{$1$} & \multicolumn{1}{c|}{$1$} & \multicolumn{1}{c|}{} \\ \hline
\multicolumn{1}{|c||}{\llap{$-$}1} & \multicolumn{1}{c|}{} & \multicolumn{1}{c|}{} & \multicolumn{1}{c|}{$I_0+3$} \\ \hline
\multicolumn{1}{|c||}{\llap{$-$}2} & \multicolumn{1}{c|}{} & \multicolumn{1}{c|}{} & \multicolumn{1}{c|}{$I_0+3$} \\ \hline
 &&& \\
\end{tabular}%
\end{minipage}\quad
%%%%%%%%%%%%%%%%%%%%%%%%%%
\begin{minipage}[t]{3.4cm}
\small
\begin{tabular}{cccc}
j=1 &  &  &  \\ \hline
\multicolumn{1}{|c||}{t\bs s} & \multicolumn{1}{c|}{0} & \multicolumn{1}{c|}{1} & \multicolumn{1}{c|}{2} \\ \hline\hline
\multicolumn{1}{|c||}{0} & \multicolumn{1}{c|}{$1$} & \multicolumn{1}{c|}{} & \multicolumn{1}{c|}{} \\ \hline
\multicolumn{1}{|c||}{\llap{$-$}1} & \multicolumn{1}{c|}{} & \multicolumn{1}{c|}{} & \multicolumn{1}{c|}{} \\ \hline
\multicolumn{1}{|c||}{\llap{$-$}2} & \multicolumn{1}{c|}{} & \multicolumn{1}{c|}{} & \multicolumn{1}{c|}{$I_0+2$} \\ \hline
\multicolumn{1}{|c||}{\llap{$-$}3} & \multicolumn{1}{c|}{} & \multicolumn{1}{c|}{} & \multicolumn{1}{c|}{$I_0+1$} \\ \hline
\end{tabular}%
\end{minipage}\quad
%%%%%%%%%%%%%%%%%%%
\begin{minipage}[t]{3.4cm}
$I_0 := (p-9)/2 \in \mathbb{Z}_{\ge 0}$
\end{minipage}%DUMMY
\caption{}\label{Tab:tss1}
\end{center}
\end{table}

When we finish writing down the possible $\KH(K)$, %
we will find that the support of possible $\KH(K)$ is in the two diagonals $\deg = 2\dim+2 \pm 1$, %
and all possible generators are shown in Table \ref{Tab:pKH1}. %
In this table, a possible generator is denoted in a bracket, like [1], [$I_0+3$]. %
So we can easily see that the support of $\KH(K)$ is also in the same two diagonals. %
In other words, $K$ is H-thin. %
From this result, we obtain $s(K)=2$ immediately. %
In this case, there is no necessity of support of Lee's theory to determine $s(K)$.

\begin{table}[htbp]
\begin{center}
\small
\begin{tabular}{|c||c|c|c|c|c|c|} \hline
deg\bs dim & \llap{$-$}1 & 0 & 1 & 2 & 3 & 4 \\ \hline\hline
9 &  &  &  &  & $I_0+3+[1]$ & $[I_0+4]$ \\ \hline
7 &  &  &  & $I_0+4+[1]$ & $[I_0+3]+1$ &  \\ \hline
5 &  &  & $I_0+2+[1]$ & $[I_0+4]+1$ &  &  \\ \hline
3 &  & $I_0+3+[1]$ & $[I_0+3]+1$ &  &  &  \\ \hline
1 & $I_0+1$ & $I_0+2+1$ &  &  &  &  \\ \hline
\llap{$-$}1 & $I_0+3$ &  &  &  &  &  \\ \hline
\end{tabular}
\end{center}
\caption{}\label{Tab:pKH1}
\end{table}

We use Lee's theory to determine $\KH(K)$. %
Because $\KH(K)$ is in the two diagonals mentioned above, Lee's theory implies that %
\begin{align}
 \mathrm{rank}\KH^{i}_{2i+1}(K) &= \mathrm{rank}\KH^{i+1}_{2i+5}(K) \text{\quad if } i \neq 0, -1 \\
 \mathrm{rank}\KH^{0}_{1}(K) -1 &= \mathrm{rank}\KH^{1}_{5}(K) \\
 \mathrm{rank}\KH^{-1}_{-1}(K) &= \mathrm{rank}\KH^{0}_{3}(K) -1 . \label{eq:+3}
\end{align}
Hence by (\ref{eq:+3}) the possible generator at bi-degree $(0,3)$ is indeed a generator. %
Then we look back Turner's $E_1 (j=3)$. %
We will find that all possible generators at bi-degree $(1,3)$ are also generators. %
In this way, all possible generators at bi-degree $(*,9)$ are indeed generators, %
and exactly one possible generator is fake at bi-degree $(1,5)$,$(2,5)$,$(2,7)$ and $(3,7)$, respectively. %
Then the homology turns out to be as shown in Table \ref{Tab:KH1}. %

\begin{table}[ht]
\begin{center}
\small
\begin{tabular}{|c||c|c|c|c|c|c|} \hline
deg\bs dim & \llap{$-$}1 & 0 & 1 & 2 & 3 & 4 \\ \hline\hline
9 &  &  &  &  & $I_0+4$ & $I_0+4$ \\ \hline
7 &  &  &  & $I_0+4$ & $I_0+3$ &  \\ \hline
5 &  &  & $I_0+2$ & $I_0+4$ &  &  \\ \hline
3 &  & $I_0+4$ & $I_0+4$ &  &  &  \\ \hline
1 & $I_0+1$ & $I_0+3$ &  &  &  &  \\ \hline
\llap{$-$}1 & $I_0+3$ &  &  &  &  &  \\ \hline
\end{tabular}
\end{center}
\caption{}\label{Tab:KH1}
\end{table}

\begin{table}[htb]
\begin{center}
\small
%$\KH(P(p,-q,-2))= \KH(P(p,-q,0))[2]\{4\} \oplus$
\begin{tabular}{|c||c|c|c|c|c|c|} \hline
deg\bs dim & $-1$ & 0 & 1 & 2 & 3 & 4 \\ \hline\hline
9 &  &  &  &  & +1 & 0 \\ \hline
7 &  &  &  & 0 & 0 &  \\ \hline
5 &  &  & 0 & 0 &  &  \\ \hline
3 &  & +1 & +1 &  &  &  \\ \hline
1 & 0 & +1 &  &  &  &  \\ \hline
\llap{$-$}1 & 0 &  &  &  &  &  \\ \hline
\end{tabular}
\end{center}
\caption{}\label{Tab:KH1inc}
\end{table} %

$\KH(K)$ is almost same as $\KH(P(p,-q,0))$ shifted by bi-degree $(2,4)$. %
Table \ref{Tab:KH1inc} shows the difference between the two homologies. %
From this viewpoint, we obtain (\ref{eq:001}). 
\end{proof}

%%%%%%%%%%%%%%%%%%%%%%%%%%%%%%%%%%%%%%%% Step 2 in case r \ge 4 %%%%%%%%%%%

\begin{prop}[Step 2 in case $r \ge 4$]\label{prop:step2ge4}
If $p:\text{odd} \ge 9$, $q = p-2 $, $r:\text{even} \ge 4$,
\begin{multline}
\Kh_{\mathbb{Q}}(K) = q^{2r} t^{r} \Kh_{\mathbb{Q}}(P(p,-q,0)) + q^1+q^3 + q^3 t^1 \\
 + \sum_{k=1}^{r/2-1} q^{4k}t^{2k}(q^1+q^3)(1+q^2 t) + q^{2r+5}t^{r+1} \label{eq:002} 
\end{multline}
\begin{equation}
s(K) = 2
\end{equation}
and $K$ is H-thin.
\end{prop}

\begin{proof}
We use Proposition \ref{Prop:Tur1} (in case $m=r$).
Let $D = P(p,-q,-r)$ and its orientation be as shown in Figure \ref{fig:ori3}, %
$D_{(k)} = P(p,-q,-(r-k)) $ and its orientation as shown in Figure \ref{fig:ori4}, %
$\tD_{(k)} = P(p,-q,-(r-k)^{\vee}_{\wedge}) $ (which inherits the orientation from $D_{(k-1)}$) for $1 \le k \le r$. %
Then $ \ta_k = 0, \tb_k = -1, A_{2l+1} = -2l-3, B_{2l+1} = -4l-8, A_{2l} = -2l, B_{2l} = -4l $ %
for $1\le k\le r, 0\le l< r/2 $, %
$A_{r} = -r, B_{r} = -2r$.

It is easy to see $ \tD_{(2l+1)} \simeq T(2,2) $ and $ \tD_{(2l)} $ is equivalent to the link $L$ obtained from $T(2,2)$ %
by changing the orientation of one component. %
We can easily get its homology from the homology of $T(2,2)$ by shifting indexes. %

The spectral sequence is
\begin{align}
E_1^{2l,t} &= \KH^{t}_{j-4l-1}(T(2,2)) & 0 \le l < r/2 \\
E_1^{2l+1,t} &= \KH^{t-2}_{j-4l-9}(L) = \KH^{t}_{j-4l-3}(T(2,2)) & 0 \le l < r/2 \\
E_1^{r,t} &= \KH^{t}_{j-2r}(P(p,-q,0)) &
\end{align}
In other words, 
\begin{equation}
E_1^{s,t} = \KH^{t}_{j-2s-1}(T(2,2)) \qquad \text{for} \quad 0 \le s \le r-1 \label{eq:015}
\end{equation}
This spectral sequence is also simple in some sense.

{\bf Claim 1}: %
Suppose that $j \in \mathbb{Z}$ fixed. %
Then the support of $E_1$ is in two diagonals $s+t=(j-2\pm 1)/2$. %
i.e. $E_1^{s,t} = 0$ if $s+t \neq (j-2 \pm 1)/2$.

{\bf Proof of Claim 1)} %
If $j$ is even, then $E_1=0$ and there is nothing to prove. %
Suppose $j$ is odd. %
Since %
\begin{equation}
\KH^a_b (T(2,2)) \neq 0
\end{equation}%
if and only if $(a,b)=(0,0),(0,2),(2,4),(2,6)$, %
we obtain
\begin{equation}
E_1^{s,t} = 0 \text{ if } (s,t) \neq (\tfrac{j-1}{2},0),(\tfrac{j-3}{2},0),(\tfrac{j-5}{2},2),(\tfrac{j-7}{2},2)
\end{equation}
for $0 \le s \le r-1$. %
For $s=r$, by Proposition \ref{prop:step1} %
we obtain that the support of $\KH(P(p,-q,0))$ is in two diagonals: $\deg = 2\dim +2 \pm 1$. %
Hence we have
$E_1^{r,t} = 0$ if $j-2r \neq 2t+2 \pm 1$.

This simplicity means the support of the possible $\KH(K)$ is also in the two diagonals $\deg = 2\dim +2 \pm 1$, %
hence $K$ is H-thin, and $s(K)=2$.

To determine $\KH(K)$, we write down the Turner's $E_1$ and possible $\KH$ explicitly, %
and use two spectral sequences alternately. %
It seems to take long time, but the well-regulatedness (\ref{eq:015}) and ``H-thinness" (Claim 1) of $E_1$ make it easier. %
We have to consider several cases: $r\le p-5,r=p-3,r=p-1,r=p+1$, and $r\ge p+3$. %
We omit details. 

\begin{table}[hbtp]
\begin{center}
\small
%$\KH(P(p,-q,-r))= \KH(P(p,-q,0))[r]\{2r\} \oplus$
\begin{tabular}{|c||c|c|c|c|c|c|c|c|} \hline
deg\bs dim & 0 & 1 & \dots & $l$ & $l+1$ & \dots & $r$ & $r+1$ \\ \hline\hline
$2r+5$ &  &  &  &  &  &  &  & +1 \\ \hline
$2r+3$ &  &  &  &  &  &  & 0 & 0 \\ \hline
 &  &  &  &  &  & $*$ & 0 &  \\ \hline
$2l+5$ &  &  &  &  & +1 & $*$ &  &  \\ \hline
$2l+3$ &  &  &  & +1 & +1 &  &  &  \\ \hline
$2l+1$ &  &  & $*$ & +1 &  &  &  &  \\ \hline
 &  & 0 & $*$ &  &  &  &  &  \\ \hline
3 & +1 & +1 &  &  &  &  &  &  \\ \hline
1 & +1 &  &  &  &  &  &  &  \\ \hline
\end{tabular}\\
$0<l<r$, $l$ : even
\end{center}
\caption{}\label{Tab:KH2inc}
\end{table} %

In all cases, %
$\KH(K)$ is almost same as $\KH(P(p,-q,0))$ shifted by bi-degree $(r,2r)$. %
Table \ref{Tab:KH2inc} shows the difference between the two homologies. %
From this viewpoint, we obtain (\ref{eq:002}).

\end{proof}

\section{Some remarks about Theorem 1.1}

\subsection{Further calculation}

\begin{rem}
There exists the case that we can't determine $\KH(K)$ by this technique. %
Even in such case, however, we may be able to determine $s(K)$. %
Example \ref{ex1} gives an example of such case.
\end{rem}

\begin{ex}\label{ex1}
If $K=P(9,-5,-2)$, %
above approach is not successful. %
But we can determine $s(K)$.
\end{ex}

One of the reasons of unsuccessfulness of calculation is %
that there is a ``double knight move" type piece of possible generators. %
There are two possibility: all four possible generators are indeed generators, or all are fake. %
We can not determine which possibility holds. Hence we need another approach. %

To determine $s(K)$ is, however, possible by this approach. %
Because possible $\KH^{0}_{j}(K)$ vanishes if $j \neq 3,5$, %
we obtain $s(K)=4$.

\begin{table}[htbp]
\begin{center}
\small
\begin{tabular}{|c|c|c|c|} \hline
 & [1] & [1] \\ \hline
 &  &  \\ \hline
[1] & [1] &  \\ \hline
\end{tabular}
\end{center}
\caption{``double knight move" type piece of possible generators}
\end{table} %

\subsection{The (sharper) slice-Bennequin type inequality for $s(K)$}\label{subsec2}

Let us focus on $s(K)$, not on $\KH(K)$. %
There are several slice-Bennequin type inequalities for $s(K)$. %
\begin{prop}[Shumakovitch\cite{Shu},Plamenevskaya\cite{Pla}]
For any knot $K$ and its diagram $D_K$, we have
\begin{equation}
s(K) \ge w(D_K) - O(D_K) + 1 \label{eq:010}
\end{equation}
where $w(D_K)$ and $O(D_K)$ are the writhe of $D_K$ and the number of Seifert circle of $D_K$, respectively.
\end{prop}
\begin{define}[Kawamura\cite{Kawamura}]
Let $D$ be an oriented diagram and $S$ a Seifert circle of $D$. %
A Seifert circle $S$ is called {\it strongly negative} if it is adjacent at least two negative crossings but no positive crossings. %
A Seifert circle $S$ is called {\it non-negative} if it is not strongly negative. 
\end{define}
\begin{prop}[\cite{Kawamura}]\label{prop:ssB}
Let $K$ be a knot and $D_K$ be a diagram of $K$. %
If $D_K$ has at least one non-negative Seifert circle, we have
\begin{equation}
s(K) \ge w(D_K) - (O_\ge(D_K)-O_<(D_K)) + 1 \label{eq:011}
\end{equation}
where $O_<(D_K)$ and $O_\ge(D_K)$ are the number of strongly negative and non-negative Seifert circle of $D_K$, respectively.
\end{prop}

In our case $K=P(p,-q,-r)$ for $p$ : odd $\ge 9$, $q=p-2$, $r$ : even $\ge 2$, %
and $D_K$ is as shown in Figure \ref{fig:ori3}, %
$D_K$ and $D_K^!$ have some non-negative Seifert circle, and %
$w(D_K) = r+2$, $O(D_K) = r+1$, $O_<(D_K) = 0$, $O_\ge(D_K) = r+1$, %
$w(D_K^!) = -r-2$, $O(D_K^!) = r+1$, $O_<(D_K^!) = r-1$, $O_\ge(D_K^!) = 2$. %
By (\ref{eq:010}), we obtain
\begin{align}
s(K) &\ge 2 \\
s(K^!) &\ge -2r-2
\end{align}
by using the property of $s(K)$ that $s(K^!)=-s(K)$, we obtain %
\begin{equation}
2 \le s(K) \le 2r+2.
\end{equation}
Similarly, by (\ref{eq:011}), we obtain
\begin{align}
s(K) &\ge 2 \\
s(K^!) &\ge -4
\end{align}
and we obtain %
\begin{equation}
2 \le s(K) \le 4.
\end{equation}

These estimations do not enable us to determine $s(K)$. %
But by the above calculation, we have obtained $s(K)=2$.

\section{Proof of Theorem 1.2 and 1.3}

Let us survey the estimation of the value of $s$-invariant of pretzel knots with three pretzels %
before proving the theorems. %
We suppose all pretzel has at least two half twists. %
By Proposition \ref{prop:property}, it is sufficient to consider the following five cases: \\
(E) Let $p,q$ be {\it positive} odd numbers and $r$ be a {\it positive} even number. \\
(E0) $K=P(p,q,r)$. 
(E1) $K=P(p,q,-r)$. 
(E2) $K=P(p,-q,-r)$. \\
(O) Let $p,q,r$ be {\it positive} odd numbers.\\
(O0) $K=P(p,q,r)$. 
(O1) $K=P(p,q,-r)$. \\

(E1) have a positive diagram, hence by \cite{Ras}, we obtain $s(K)=p+q$. %
Similarly, (O0) have a negative diagram, hence $K^!$ have a positive one, so we obtain $s(K)=-2$. %
In the other three cases (E0), (E2) and (O1), Proposition \ref{prop:ssB} gives us estimations of $s(K)$, %
as in subsection \ref{subsec2}. %
For details, see Table \ref{tab:estsE} and Table \ref{tab:estsO}. %
These tables show alternatingness of $K$, positiveness (or negativeness) of $K$, %
and the values of $s(K)$ and the signature $\sigma(K)$ of $K$. %
If $K$ is alternating, then $s(K)=-\sigma(K)$.
In these table, we set $\omega\big(P(i,j,k)\big) := ij+jk+ki$.

\begin{table}[hbtp]
\begin{center}
\footnotesize
\begin{tabular}{c|ccccc} \hline
$p,q$ : odd $\ge 3$\\$r$ : even $\ge 2$ & alternating & posi/nega & $s(K)$ & $\sigma(K)$ & \\ \hline
$P(p,q,r)$ & alternating & ? & $p+q-2$ & $-(p+q-2)$ & \\ \hline
$P(p,q,-r)$ & ? & positive & $p+q$ & $-p-q$ & if $\omega<0$ \\
            & non-alt &     &       & $-p-q+2$ & if $\omega>0$ \\ \hline
$P(p,-q,-r)$ & ? & ? & $p-q-2\le s(K) \le p-q$ & $-(p-q)$ & if $\omega<0$ \\
            &   &   &              & $-(p-q-2)$ & if $\omega>0$ \\ \hline
\end{tabular}
\end{center}
\caption{}\label{tab:estsE}
\end{table} %
\begin{table}[hbtp]
\begin{center}
\footnotesize
\begin{tabular}{c|ccccc} \hline
$p,q,r$ : odd $\ge 3$ & alternating & posi/nega & $s(K)$ & $\sigma(K)$ & \\ \hline
$P(p,q,r)$ & alternating & negative & -2 & 2 & \\ \hline
$P(p,q,-r)$ & ? & ? & $-2\le s(K) \le 0$ & 0 & if $\omega<0$ \\
            &   &   &                    & 2 & if $\omega>0$ \\ \hline
\end{tabular}
\end{center}
\caption{}\label{tab:estsO}
\end{table} %

In case (E0), $s(K)$ has been already determined. %
In case (E2), we determine $s(K)$ in Theorem \ref{Thm:main2}. %
In some case in (E3), we determine $s(K)$ in Theorem \ref{Thm:main3}. %
The next lemma plays an important part as computational shortcut.

\begin{lem}\label{lem1}
If $p,q$ are odd numbers $\ge 3$, then $K=T(2,p)\#T(2,-q)$ is H-thin.
\end{lem}

\begin{proof}
We use Turner's spectral sequence with following settings. (Same as in the proof of Proposition \ref{prop:step1}.) %
Let $D = D_{(0)} = \tD_{(0)} = P(p,-q,0)$ and its orientation be as shown in Figure \ref{fig:ori1}. %
For $1 \le k \le q-1$, $D_{(k)} := P(p,-(q-k),0)$ (which inherits the orientation from $D_{(k-1)}$), %
$\tD_{(k)} := P(p,-(q-k)^{\vee}_{\wedge},0)$ and its orientation be as shown in Figure \ref{fig:ori2}. %
It is easy to see that $\tD_{(k)} \simeq T(2,p)$ for $1 \le k \le q-1$, $D_{(q-1)} \simeq T(2,p)$, %
$\ta_s = q-s+1$, $\tb_s = 3q-3s+2$, $A_s = 0$, $B_s = s$ for $1 \le s \le q-1$.

The spectral sequence is
\begin{align}
E_1^{s,t} &= \KH^{t+q}_{j+3q-2s-1}(T(2,p)) \text{ if } s=0,\dots,q-2\\
E_1^{q-1,t} &= \KH^{t+q-1}_{j+q-1}(T(2,p)).
\end{align}

Here we note that for $p \ge 2$, support of the homology of $(2,p)$-torus knot is given by
\begin{equation}
\supp\KH(T(2,p)) \subset \bigl\{ \deg = p-1 \pm 1 +2\dim, 0 \le \dim \le p \bigr\}. \label{eq:030}
\end{equation}
Hence the support of the spectral sequence is given by %
\begin{align}
\supp E_1^{s,t} &\subset \Bigl\{ t = \frac{j+q-p\pm1}{2} -s, -q \le t \le p-q \Bigr\} \text{ if } s=0,\dots,q-2\\
\supp E_1^{q-1,t} &\subset \Bigl\{ t= \frac{j-q-p\pm1}{2} -1, -q+1 \le t \le p-q+1 \Bigr\}.
\end{align}
In other words, $E_1^{s,t}=0$ if $s+t \neq \frac{j+q-p\pm1}{2}$. %
From this we conclude that $\KH^i_j(K) = 0$ if $j \neq 2i +p-q \pm1 $. %
This completes the proof.
\end{proof}

The following proposition partially proves Theorem \ref{Thm:main2}. %
Note that there is an additional condition $p>q$.
\begin{prop}
If $p,q$ are odd numbers $\ge 3$, $p>q$, $r$ is a positive even number and $K=P(p,-q,-r)$, then %
\begin{equation}
s(K)=p-q.
\end{equation}
\end{prop}

\begin{proof}
We use Turner's spectral sequence with following settings. (Same as in the proof of Proposition \ref{prop:step2ge4}.) %
Let $D = P(p,-q,-r)$ and its orientation be as shown in Figure \ref{fig:ori3}, %
$D_{(k)} = P(p,-q,-(r-k)) $ and its orientation as shown in Figure \ref{fig:ori4}, %
$\tD_{(k)} = P(p,-q,-(r-k)^{\vee}_{\wedge}) $ (which inherits the orientation from $D_{(k-1)}$) for $1 \le k \le r$. %
Then $ \ta_k = 0$, $\tb_k = -1$, $A_{2l+1} = -2l-1-p+q$, $B_{2l+1} = -4l-2-3p+3q$, $A_{2l} = -2l$, $B_{2l} = -4l $ %
for $1\le k\le r, 0\le l< r/2 $, %
$A_{r} = -r$, $B_{r} = -2r$.

It is easy to see $ \tD_{(2l+1)} \simeq T(2,p-q) $ and $ \tD_{(2l)} $ is equivalent to the link $L$ obtained from $T(2,p-q)$ %
by changing the orientation of one component. %
We can easily get its homology from the homology of $T(2,p-q)$ by shifting indexes: %
\begin{equation}
\KH^{i}_{j}(L) \simeq \KH^{i+(p-q)}_{j+3(p-q)}(T(2,p-q))
\end{equation}

The spectral sequence is
\begin{align}
E_1^{s,t} &= \KH^{t}_{j-2s-1}(T(2,p-q)) \qquad \text{for} \quad 0 \le s \le r-1 \\
E_1^{r,t} &= \KH^{t}_{j-2r}(P(p,-q,0)).
\end{align}
From (\ref{eq:030}), %
we obtain that $ E_1^{s,t}|_{s+t=0,0\le s \le r-1} \neq 0 $ only if $s=0$, $j=p-q\pm1$. %
Moreover, from Lemma \ref{lem1}, %
we obtain $ E_1^{r,t}|_{r+t=0} \neq 0 $ only if $ j=p-q\pm1 $. %
Hence we conclude $ E_1^{s,t}|_{s+t=0} \neq 0 $ only if $j=p-q\pm1$, %
which proves the proposition.
\end{proof}

To complete the proof of Theorem \ref{Thm:main2}, %
we use the following fact. %

\begin{lem}[Livingston\cite{Liv}, J. Rasmussen\cite{Ras}]\label{lemLR}
Suppose $K_{+}$ and $K_{-}$ are knots that differ by a single crossing change %
--- from a positive crossing in $K_{+}$ to a negative one in $K_{-}$. Then 
\[ s(K_{-}) \le s(K_{+}) \le s(K_{-})+2. \]
\end{lem}

The rest of the theorem will be proved by induction. %
Let $K_{+}=P(p,-q,-r)$, $K_{-}=P(p-2,-q,-r)$. %
Suppose that $s(P(p,-q,-r))=p-q$. %
From the lemma above, we obtain %
\begin{equation}
 p-q-2 \le s(K_{-}) \le p-q. 
\end{equation}
On the other hand, from the slice-Bennequin type inequality for $s(K)$ (see Table \ref{tab:estsE}), we obtain %
\begin{equation}
 p-q-4 \le s(K_{-}) \le p-q-2. 
\end{equation}
Then we obtain $s(K_{-})=p-q-2$, %
and we have completed the proof of Theorem \ref {Thm:main2}.

Let us prove Theorem \ref{Thm:main3}. %
We begin by proving in the case $p=q+2$ and $r$ is arbitrary. %

\begin{prop}
If $p,q$ are odd numbers $\ge 3$, $p=q+2$, $r$ is a positive odd number and $K=P(p,-q,-r)$, then %
\begin{equation}
s(K)=0.
\end{equation}
\end{prop}

\begin{proof}
We use the spectral sequence with following settings. %
Let $D=P(p,-q,-r)$, $D_{(1)}=P(p,-q,-(r-1)$ and its orientation be as shown in Figure \ref{fig:ori4} (set $k=0,1$), %
$\tD=(p,-q,-(r-1)^{\vee}_{\wedge})$ (which inherits the orientation from $D_{(0)}=D$). %
We note that $r$ is an odd number. %
Then $ \ta_1=0$, $\tb_1=-1$, $A_1=1$, $B_1=4$. %

In this case, our spectral sequence is %
\begin{align}
E_1^{0,t} &= \KH^{t}_{j-1}(T(2,-2)) =\KH^{t+2}_{j+5}(T(2,2)) \\
E_1^{1,t} &= \KH^{t+2}_{j+4}(P(p,-q,-(r-1))).
\end{align}
It is sufficient to focus on $ E_1^{s,t}|_{s+t=0} $. %
It is easy to check that $ E_1^{s,t}|_{s+t=0} = 0 $ if $j \neq \pm1$, %
which completes the proof. %
If $s=0$, there is nothing to comment. %
If $s=1$, it follows from Theorem \ref{Thm:main}; %
we obtain $\KH^{1}_{j+4}(P(p,-q,-(r-1))=0$ if $j \neq \pm 1$ %
from Theorem \ref{Thm:main}.
\end{proof}

By using Lemma \ref{lemLR} inductively, the condition $p=q+2$ is weakened to $p>q$. %
Let us see the detail. %
Let $K_{+}=P(p,-q,-r)$, $K_{-}=P(p+2,-q,-r)$. %
Suppose that $s(P(p,-q,-r))=0$.
Lemma \ref{lemLR} gives us 
$ -2 \le s(K_{-}) \le 0$. %
On the other hand, the slice-Bennequin type inequality for $s(K)$ gives us 
$ 0 \le s(K_{-}) \le 2$. %
Then we obtain $s(K_{-})=0$. %
Moreover, the condition $p>q$ is weakened to $p>\min\{q,r\}$ %
because of Proposition \ref{prop:property}. %

It remains to consider the case $p \le \min\{q,r\}$. %
We give the partial result for the case $p < \min\{q,r\}$, %
which completes the proof of Theorem \ref{Thm:main3}.

\begin{prop}
If $p,q,r$ are odd numbers $\ge 3$, $p<\min\{q,r\}$, then %
\begin{equation}
s(K)=2.
\end{equation}
\end{prop}

\begin{proof}
Let us suppose $p<q$ and %
consider Turner's spectral sequence with following settings. %
Let $D_{(k)} = P(p,-q,-(r-k)) $ and its orientation as shown in Figure \ref{fig:ori4}, for $0 \le k \le r$. %
Let $\tD_{(k)} = P(p,-q,-(r-k)^{\vee}_{\wedge}) $ (which inherits the orientation from $D_{(k-1)}$) for $1 \le k \le r$. %
Then $ \ta_k = 0$, $\tb_k = -1$, $A_{2l+1} = -2l-1+p-q$, $B_{2l+1} = -4l-2+3p-3q$, $A_{2l} = -2l$, $B_{2l} = -4l $ %
for $1\le k\le r, 0\le l< r/2 $. %
Especially, $A_{r} = -r+p-q$, $B_{r} = -2r+3(p-q)$. %
It is easy to see that $\tD_{(k)} \simeq T(2,p-q)$.

Our spectral sequence is
\begin{align}
E_1^{s,t} &= \KH^{t+p-q}_{j-4l-1+3(p-q)}(T(2,p-q)) \text{ if } s=0,\dots,r-1\\
E_1^{r,t} &= \KH^{t+p-q}_{j-2r+3(p-q)}(T(2,p)\#T(2,-q)).
\end{align}
We focus on $ E_1^{s,t}|_{s+t=0} $. %
$E_1^{0,0} = 0$ if $j \neq 1,3 $. %
Since $\KH^i(T(2,p-q))=0$ if $i<p-q$, %
$E_1^{s,-s} = 0$ for $1 \le s \le r-1$ and all $j$. %
Here we suppose $p<r$, then $E_1^{r,-r} = 0$ for all $j$ because $\KH^i(T(2,p)\#T(2,-q))=0$ if $i<-q$. %
Hence we find that $\KH^{0}_{j}(D)=0$ if $j \neq 1,3$. %
So we conclude that $s(P(p,-q,-r))=2$ if $p<\min\{q,r\}$.
\end{proof}

\begin{rem}
We have not yet determined $s(K)$ in case $ p=\min\{q,r\} $.
\end{rem}

\begin{rem}
By Lemma \ref{lemLR}, the following equations are hold:
if $p,q,r$ are odd numbers, then
\begin{align}
P(p+2,-q,-r) &= P(p,-q,-r) -0 \text{ or } -2 \\
P(p,-(q+2),-r) &= P(p,-q,-r) +0 \text{ or } +2 \\
P(p,-q,-(r+2)) &= P(p,-q,-r) +0 \text{ or } +2.
\end{align}
These inequalities played important role in the proof of Theorem \ref{Thm:main3}, %
but they do not help to weaken the assumption of Theorem \ref{Thm:main3} anymore.
\end{rem}

{\bf Acknowledgement.}
The author gratefully acknowledges the many helpful suggestions of Professor Kohno during the preparation of the paper.

%%%%%%%%%%%%%%%%%%%%%%%%%%%%%%%%%%%%%%%%%%%%%%%%%%%%%%%%%%%%%%%%%%%%%%%%

{\footnotesize
{\sc Graduate School of Mathematical Sciences, University of Tokyo, 3-8-1 Komaba Meguro, Tokyo 153-8914, Japan}

E-mail address : \texttt{rsuzuki@ms.u-tokyo.ac.jp}
}

%%%%%%%%%%%%%%%%%%%%% huge figures : Table 1 %%%%%%%%%%%%%%%%%%%%%%%%%

\begin{landscape}

\begin{table}[htbp]
{\small $
\begin{array}{c||c|c|c|c|c|c|c|c|c|c|c|c|c|c|c|}
\rm{deg\bs dim}&\ldots&-2&-1& 0 & 1 & 2 &\ldots& m &m+1&m+2&\ldots&p-1& p \\ \hline \hline
2p+3         &\phantom{I_m+0}&  &  &   &   &   &      &   &   &   &      &   & 1 \\ \hline
2p+1         &      &  &  &   &   &   &      &   &   &   &      &   &   \\ \hline
\vdots       &      &  &  &   &   &   &      &   &   &   &\bdots&   &   \\ \hline
5+2m         &      &  &  &   &   &   &      &   &I_m+2&I_m+1&      &   &   \\ \hline
3+2m         &      &  &  &   &   &   &      &I_m+1&I_m&   &      &   &   \\ \hline
\vdots       &      &  &  &   &   &   &\bdots&   &   &   &      &   &   \\ \hline
5            &      &  &  &   &I_0+3&I_0+4&      &   &   &   &      &   &   \\ \hline
3            &      &  &  &I_0+4&I_0+3&   &      &   &   &   &      &   &   \\ \hline
1            &      &  &I_0+2&I_0+4&   &   &      &   &   &   &      &   &   \\ \hline
-1           &      &I_0+3&I_0+3&   &   &   &      &   &   &   &      &   &   \\ \hline
\vdots       &\bdots_*&&  &   &   &   &\phantom{I_m+0}&   &   &   &\phantom{I_m+0}&\phantom{I_m+0}&\phantom{I_m+0}\\ \hline
\end{array} $\\
 \phantom{M}

$
\begin{array}{c||c|c|c|c|c|c|c|c|c|c|c|c|c|c|c|}
\rm{deg\bs dim} &-p+2&-p+3&-p+4&-p+5&\ldots&n-2&n-1& n &\ldots&-2&-1& 0 \\ \hline \hline
-1           &    &    &    &    &      &   &   &   &      &I_0+3&I_0+3&   \\ \hline
\vdots       &    &    &    &    &      &   &   &   &\bdots_*&&  &   \\ \hline
1+2n         &    &    &    &    &      &   &I_n&I_n+1&      &  &  &   \\ \hline
-1+2n        &    &    &    &    &      &I_n+1&I_n+2&   &      &  &  &   \\ \hline
\vdots       &    &    &    &    &\bdots&   &   &   &      &  &  &   \\ \hline
-2p+11       &    &    &    &  1 &      &   &   &   &      &  &  &   \\ \hline
-2p+9        &    &  1 &  2 &    &      &   &   &   &      &  &  &   \\ \hline
-2p+7        &    &    &    &    &      &   &   &   &      &  &  &   \\ \hline
-2p+5        &  1 &    &    &    &\phantom{I_m+0}&   &   &   &\phantom{I_m+0}&  &  &\phantom{I_m+0}\\ \hline
\end{array}
$}%small

$ m : \rm{even} , 2 \le m \le p-3 $, 
$ n : \rm{even} , -p+7 \le n \le -2 $, 
and $ I_m = \frac{p-m-3}{2}, I_0 = \frac{p-9}{2}, I_n = \frac{p+n-5}{2} $. 

\caption{The rational Khovanov homology of $P(p,-q,0)$ for $p:\text{odd} \ge 7, q = p-2 $}\label{TabHuge1}

\end{table}

\end{landscape}

%%%%%%%%%%%%%%%%%%%%% huge figures : Table 2 %%%%%%%%%%%%%%%%%%%%%%%%%
\begin{landscape}

\begin{table}[htbp]
{\footnotesize $
\begin{array}{c||c|c|c|c|c|c|c|c|c|c|c|c|c|c|c|}
\rm{deg\bs dim} &\ldots&-2&-1& 0 & 1 & 2 &\ldots& m &m+1&m+2&\ldots&p-1& p \\ \hline \hline
2p+3         &      &  &  &   &   &   &      &   &   &   &      &   & 1 \\ \hline
2p+1         &      &  &  &   &   &   &      &   &   &   &      &   &   \\ \hline
\vdots       &      &  &  &   &   &   &      &   &   &   &\bdots&   &   \\ \hline
5+2m         &      &  &  &   &   &   &      &   &[I_m+1]+1&[I_m+1]&      &   &   \\ \hline
3+2m         &      &  &  &   &   &   &      &[I_m]+1&[I_m]&   &      &   &   \\ \hline
\vdots       &      &  &  &   &   &   &\bdots&   &   &   &      &   &   \\ \hline
5            &      &  &  &   &[I_0+3]&[I_0+3]+1&      &   &   &   &      &   &   \\ \hline
3            &      &  &  &1+[I_0+3]&[I_0+3]&   &      &   &   &   &      &   &   \\ \hline
1            &      &  &[I_0+2]&[I_0+2]+2&   &   &      &   &   &   &      &   &   \\ \hline
-1           &      &[I_0+2]+1&[I_0+2]+1&   &   &   &      &   &   &   &      &   &   \\ \hline
\vdots       &\bdots&  &  &   &   &   &      &   &   &   &      &   &   \\ \hline
\end{array} $\\
 \phantom{M}

$
\begin{array}{c||c|c|c|c|c|c|c|c|c|c|c|c|c|c|c|}
\rm{deg\bs dim} &-p+2&-p+3&-p+4&-p+5&\ldots&n-2&n-1& n &\ldots&-2&-1& 0 \\ \hline \hline
-1           &    &    &    &    &      &   &   &   &      &[I_0+2]+1&[I_0+2]+1&   \\ \hline
\vdots       &    &    &    &    &      &   &   &   &\bdots&  &  &   \\ \hline
1+2n         &    &    &    &    &      &   &[I_n+1]&[I_n+1]+1&      &  &  &   \\ \hline
-1+2n        &    &    &    &    &      &[I_n+1]&[I_n+1]+1&   &      &  &  &   \\ \hline
\vdots       &    &    &    &    &\bdots&   &   &   &      &  &  &   \\ \hline
-2p+11       &    &    & [1]&[1]+1&      &   &   &   &      &  &  &   \\ \hline
-2p+9        &    & [1]&[1]+1&    &      &   &   &   &      &  &  &   \\ \hline
-2p+7        & [1]& [1]&    &    &      &   &   &   &      &  &  &   \\ \hline
-2p+5        &  1 &    &    &    &      &   &   &   &      &  &  &   \\ \hline
\end{array}
$}%small

$ m : \rm{even} , 2 \le m \le p-3 $, 
$ n : \rm{even} , -p+7 \le n \le -2 $, 
and $ I_m = \frac{p-m-3}{2}, I_0 = \frac{p-9}{2}, I_n = \frac{p+n-5}{2} $. 

A possible generator is in a bracket, like [1], [$I_0+3$].

\caption{The possible situation of the rational Khovanov homology of $P(p,-q,0)$ for $p:\text{odd} \ge 7, q = p-2 $}\label{TabHuge2}

\end{table}

\end{landscape}
%%%%%%%%%%%%%%%%%%%%% end of huge figures %%%%%%%%%%%%%%%%%%%%%%%%%

\end{document}